\documentclass[12pt]{amsart}
\usepackage[dvips]{graphicx}
\usepackage{amsmath,graphics}
\usepackage{amsfonts,amssymb}
\theoremstyle{plain}
\newtheorem{theorem*}{Theorem}
\newtheorem*{lemma*} {Lemma}
\newtheorem*{corollary*} {Corollary}
\newtheorem*{proposition*}{Proposition}
\newtheorem{conjecture*}{Conjecture}
\newtheorem{theorem}{Theorem}[section]
\newtheorem{lemma}[theorem]{Lemma}

\newtheorem{proposition}[theorem]{Proposition}

\theoremstyle{remark}
\newtheorem*{remark}{Remark}

\theoremstyle{definition}
\newtheorem{defn}[theorem]{Definition}

\def\gl{\mbox{GL}}  \def\F{\Bbb{F}} \def\Z{\Bbb{Z}}  
    
 \def\a{\alpha}   \def\bp{\begin{pmatrix}}
\def\sm{\setminus} \def\ep{\end{pmatrix}} \def\bn{\begin{enumerate}} 
   \def\en{\end{enumerate}}
\def\ba{\begin{array}} \def\ea{\end{array}}  
   \def\a{\alpha}  \def\ti{\tilde}
  \def\im{\mbox{Im}} 
  
\def\be{\begin{equation}} \def\ee{\end{equation}}

 \def\aut{\mbox{Aut}}

\def\f13{\F_{13}}

 \def\ftht{\F_{13}[t^{\pm 1}]}
 \def\deltaks5{\Delta_{K}^{5}}

    \def\f{\F} \def\GL{\mbox{GL}}   
\def\zt{\Z[t^{\pm 1}]}    
   
    \def\fr12{\frac{1}{2}} \def\z12{\Z[\fr12]} 
 \def\rf{R[F]}

\def\tpm {[t^{\pm 1}]}

\begin{document}

\title{Nontrivial Alexander polynomials of knots and links}
\author{Stefan Friedl}
\address{Universit\'e du Qu\`ebec \`a Montr\'eal, Montr\'eal, Qu\`ebec}
\email{sfriedl@gmail.com}
\author{Stefano Vidussi}
\address{Department of Mathematics, University of California,
Riverside, CA 92521, USA} \email{svidussi@math.ucr.edu} \thanks{S. Vidussi was partially
supported by NSF grant \#0629956.}
\date{June 22, 2006}
\begin{abstract}
In this paper we present a sequence of link invariants, defined from twisted Alexander
polynomials, and discuss their effectiveness in distinguishing knots. In particular, we recast and
extend by geometric means a recent result of Silver and Williams on the nontriviality of twisted
Alexander polynomials for nontrivial knots. Furthermore building on results in \cite{FV06b} we
prove that these invariants decide if a genus one knot is fibered. Finally we also show that these
invariants distinguish all mutants with up to 12 crossings.

\end{abstract} \maketitle

\section{Definition of the invariant and main results}
Let $L \subset S^3$ be an oriented $m$--component link, and denote by $X(L) = S^{3} \setminus
\nu L$ its exterior. Let $R=\Z$ or $R=\F_p$: given a representation $\a : \pi_1(X(L)) \to
\GL(R,k)$ we can consider the associated multivariable twisted Alexander polynomial
$\Delta^\a_L\in R[t_1^{\pm 1},\dots,t_m^{\pm 1}]$ (where $t_1,\dots,t_m$ correspond to a basis
of $H_1(X(L))$ determined by the meridians to each link component), well--defined up to units.
In Section \ref{sectionufd} we recall the details of the definition.

Let now $\a:\pi_1(X)\to S_k$ be a homomorphism into the symmetric group. Using the action of
$S_k$ on $R^k$ by permutation of the coordinates, we get a representation $\pi_1(X)\to S_k\to
\gl(R,k)$ that we will denote by $\a$ as well. Consider now the set of representations of $\pi_1(X(L))$ in the symmetric group, modulo conjugation:

\[ \mathcal{R}_k(L) = \{ \a: \pi_1(X(L)) \to S_k \}/\sim \]
where two representations are equivalent if they are the same up to conjugation by an element in
$S_k$. Given $\a: \pi_1(X(L)) \to S_k$ the polynomial $\Delta^\a_L$ depends only on
the equivalence class $[\a]$ of $\a$ in $\mathcal{R}_{k}(L)$. We now define  the invariant
\[ \Delta^k_L = \prod_{[\a] \in \mathcal{R}_k(L)} \Delta^\a_L \in R[t_1^{\pm 1},\dots,t_m^{\pm 1}]. \]

We will illustrate the effectiveness of this invariant by discussing some of the topological
information that it carries, and by using explicit calculations we show its ability to distinguish
many examples of inequivalent mutant knots.

Our first result relates the link invariants $\Delta_L^k$ with epimorphisms of the link group onto finite groups, which will lead to a useful
topological interpretation (cf. Lemma \ref{lemmaalexg}). Precisely, consider an epimorphism $\gamma: \pi_1(X(L)) \to G$, where $G$ is a finite group of order $k =
|G|$. Using the left action of $G$ on its group ring we can define a representation, denoted
with the same symbol, $\gamma: \pi_1(X(L)) \to \aut_{R}(R[G])\cong \mbox{GL}(R,k)$. We have the
following:

\begin{proposition} \label{prop:deltagdivides}
Let $\gamma: \pi_1(X(L)) \to G$ be a homomorphism to a finite group $G$ of order $k$. Then $\Delta_L^\gamma$ divides
$\Delta_L^{k}$.
\end{proposition}

This
relation is crucial in proving the following theorem, which shows that the sequence $\Delta_L^k$
detects the unknot and the Hopf link.

\begin{theorem}\label{mainthm}
Let $L \subset S^3$ be an oriented link which is neither the unknot nor the Hopf link (with
either orientation). Then there exists a $k$ such that $\Delta^k_L \ne \pm 1\in \Z[t_1^{\pm
1},\dots,t_m^{\pm 1}]$.
\end{theorem}

In fact we will show that if $L$ is neither the unknot nor the Hopf link,
then there exists an epimorphism $\gamma: \pi_1(X(L)) \to G$ to a finite group such that  $\Delta_L^\gamma \neq 1$.
(This result is nontrivial when $m=1$ or $2$.)
For the case of knots this provides a different approach to a recent result by Silver and
Williams \cite{SW06}.

The proof is based on the relation between twisted Alexander polynomials and covers
of the link exterior, using ideas from previous papers by the authors \cite{FV06,FV06b}, combined
with information on the topology of those covers arising from the work in \cite{Lu88,Ko87,CLR97}.

If $K$ is a fibered knot, its ordinary Alexander polynomial is monic. The following result,
combining results from \cite{FK06} and \cite{FV06b}, generalizes that assertion to $\Delta^k_K$
and shows that, at least in some cases, the converse holds true.
\begin{theorem}\label{thm:fibered}
Let $K \subset S^3$ be a fibered knot, then  $\Delta^k_K \in \Z[t^{\pm 1}]$ is monic for any $k$. Conversely, if
 $\Delta^k_K$ is monic for all $k$ and if $K$ is a genus one knot, then $K$ is
fibered.
\end{theorem}

Note that the converse also holds for knots whose exteriors has fundamental group that satisfies
suitable subgroup separability
properties. We refer the interested reader to \cite{FV06b} for details (the results in
\cite{FV06b} are only stated for closed 3--manifolds, but they also hold for 3--manifolds with
toroidal boundary).

For a knot $K$ the calculation of $\Delta_K^k$ can be done using the program \emph{KnotTwister}
\cite{Fr06}. Our computations in Section \ref{section:examples} confirm that $\Delta_K^k$ are very
strong knot invariants. For example computing $\Delta_K^{5}\in \F_{13}\tpm$ distinguishes all
pairs and triples of mutants with up to 12 crossings (cf. Section \ref{section:examples} for the
definition of mutants). In Section \ref{section:examples} we also show that $\Delta_K^4$ is not
determined by either HOMFLY polynomial, Khovanov homology or Knot Floer homology.
\\

The paper is organized as follows. In Section \ref{section:def} we give a precise definition of
twisted Alexander polynomials and discuss some basic properties. In particular we give a proof of
Proposition \ref{prop:deltagdivides}. In Section \ref{section:proofs} we give the proofs of
Theorems \ref{mainthm} and \ref{thm:fibered}. We conclude the paper in Section
\ref{section:examples} with several examples and questions.
\\

\noindent {\bf Acknowledgment.} The authors would like to thank Xiao-Song Lin for useful
discussions and Alexander Stoimenow for providing the braid descriptions for the mutants and for
giving  very helpful feedback to the program \textit{KnotTwister}.

%===================================================
\section{Twisted Alexander polynomials and finite covers} \label{section:def}
%===================================================
\subsection{Twisted Alexander modules and their polynomials} \label{sectionufd}
In this section we give the precise definition of the (twisted) Alexander polynomials. Twisted
Alexander polynomials were introduced, for the case of knots, in 1990 by Lin \cite{Li01}, and
further generalized to links by Wada \cite{Wa94}. We follow the approach taken by Cha
\cite{Ch03} and \cite{FV06}.

 For the
remainder of this section let $N$ be a  3--manifold (by which we always mean a compact, connected
and oriented 3--manifold) and denote by $H := H_1(N)/\mbox{Tor$H_1(N)$}$ the maximal free abelian quotient
of $\pi_{1}(N)$. Furthermore let $F$ be a free abelian group and let $R$ be $\Z$ or the field
$\F_p:=\Z/p\Z$ where $p$ is a prime number.

Now let $\phi \in \mbox{Hom}(H,F)$ be a non--trivial homomorphism. Through the homomorphism $\phi$,
$\pi_{1}(N)$ acts on $F$ by translations. Furthermore let $\a:\pi_1(N)\to \gl(R,k)$ be a
representation. We write $R^k[F]=R^k\otimes_R R[F]$. We get a representation
\[ \ba{rrcl}  \a\otimes \phi:&\pi_1(N)&\to& \aut(R^k[F])\\
&g&\mapsto& (\sum_i a_i\otimes f_i \mapsto \sum_i\a(g)(a_i)\otimes (f_i+\phi(g)).\ea \] We can
therefore view $R^k[F]$ as a left $\Z[\pi_1(N)]$--module. Note that this module structure
commutes with the natural  $R[F]$--multiplication on $R^k[F]$.

Let $\ti{N}$ be the universal cover of $N$. Note that $\pi_{1}(N)$ acts on the left on $\ti{N}$ as
the group of deck transformation. The chain groups $C_*(\ti{N})$ are in a natural way right
$\Z[\pi_1(N)]$--modules, with the right action on $C_{*}(\ti{N})$ defined via $\sigma \cdot g :=
g^{-1}\sigma$, for $\sigma \in C_{*}(\ti{N})$. We can form by tensor product the chain complex
$C_*(\ti{N})\otimes_{\Z[\pi_1(N)]}R^k[F]$. Now define $H_{1}(N;R^k[F]):=
H_1(C_*(\ti{N})\otimes_{\Z[\pi_1(N)]}R^k[F])$, which inherits the structure of $\rf$--module.

The $\rf$--module $H_1(N;R^k[F])$ is a finitely presented and finitely related $\rf$--module since
$\rf$ is Noetherian. Therefore $H_1(N;R^k[F])$ has a free $\rf$--resolution
\[ \rf^r \xrightarrow{S} \rf^s \to H_1(N;R^k[F]) \to 0 \]
of finite $\rf$--modules, where we can always assume that $r \geq s$.

\begin{defn} \label{def:alex} The \emph{twisted Alexander polynomial} of $(N,\a,\phi)$ is defined
to be  the order of the $\rf$--module $H_1(N;R^k[F])$, i.e. the greatest common divisor of the
$s\times s$ minors of the $s\times r$--matrix $S$. It is denoted by $\Delta_{N,\phi}^{\a}\in
\rf$.
\end{defn}

Note that this definition only makes sense since $\rf$ is a UFD. It is well--known that
$\Delta_{N,\phi}^{\a}$ is well--defined only up to multiplication by a unit in $\rf$ and its
definition is independent of the choice of the resolution.

 When $\phi$ is the identity map on $H$, we will simply write $\Delta_{N}^{\a}$. Also, we will write
${\Delta}_{N,\phi}$ in the case that $\a:\pi_1(N)\to \gl(\Z,1)$ is the trivial representation.

If $N=X(L)$ is the exterior of an oriented ordered link $L=L_1\cup\dots\cup L_m$, then we write $\Delta_{L}^\a$ for the twisted Alexander polynomial of $X(L)$. Also, we can
identify $H$ with the free abelian group generated by $t_1,\dots,t_m$ and we can view the
corresponding twisted Alexander polynomial $\Delta^\a_L$ as an element in $R[t_1^{\pm
1},\dots,t_m^{\pm 1}]$.
%\begin{remark} Note that
%we are defining the twisted Alexander polynomial as an invariant of the exterior of the link, in
%particular we are not fixing a canonical system of generators for $H_1(X(L))$, so that links with
%homeomorphic exterior have the same polynomial up to change of generators. This however, bears no
%effect on our result. \end{remark}

%===================================================
\subsection{Twisted Alexander polynomials and homomorphisms to finite groups}
Let $N$ be a  3--manifold and let $\gamma:\pi_1(N) \to G$ be an epimorphism onto a finite group
$G$ of order $k = |G|$. We get the induced regular representation $\gamma:\pi_1(N)\to G\to
\aut(R[G])$ where $g\in G$ acts on $R[G]$ by left multiplication. Since $R[G]\cong R^{|G|}$ we
can identify $\aut_{R}(R[G]) = \GL(R,k)$. It is easy to see that the isomorphism type of the
$R[H]$--module $H_1(N;R^k[H])$ does not depend on the identification $\aut_{R}(R[G]) = \GL(R,k)$.

The following lemma clearly implies Proposition \ref{prop:deltagdivides}.

\begin{lemma}\label{lemma:deltagdivides}
Let $\gamma:\pi_1(N)\to G$ be an epimorphism onto a finite group $G$ of order $k$. Then there exists a homomorphism
$\a:\pi_1(N)\to S_{k}$ such that the corresponding representation
\[ \a: \pi_1(N)\to S_{k}\to \gl(R,k) \]
is given by the regular representation $\gamma: \pi_1(N)\to G\to \gl(R,k)$.
\end{lemma}

\begin{proof}
Denote the elements of $G$ by
$g_1,\dots,g_{k}$. Since $\gamma$ defines an action on the set $G=\{g_1,\dots,g_k\}$ via left
multiplication we get an induced map $\a:\pi_1(N)\to S_{k}$. Clearly the corresponding
representation
\[ \a: \pi_1(N)\to S_k\to \gl(R,k) \]
is isomorphic to the regular representation $\gamma: \pi_1(N)\to G\to \gl(R,k)$.
\end{proof}

%===================================================
\subsection{Twisted Alexander polynomials and finite covers}
For the remainder of this section let $\gamma:\pi_1(N) \to G$ be an epimorphism onto a finite
group $G$ of order $k$, and take $R = \Z$. Denote the induced $G$--cover of $N$ by $\pi : N_{G}
\rightarrow N$. Also, denote by $H_{G}$ the maximal free abelian quotient of $\pi_{1}(N_{G})$: the map
$\pi_{*} : H_G \to H$ is easily seen to have maximal rank, hence in particular $b_1(N_G)\geq
b_1(N)$. Given any homomorphism $\phi: H \rightarrow F$ to a free abelian group $F$ we can
consider the induced homomorphism $\phi_{G} := \pi^{*}\phi : H_G \rightarrow F$. In particular,
when $\phi$ is the identity map on $H$, we have $\phi_{G} = \pi_{*} : H_G \rightarrow H$.

We can now formulate the relationship between the twisted Alexander polynomials of $N$ and
the untwisted Alexander polynomial of $N_G$.

\begin{lemma} \label{lemmaalexg}\cite{FV06}
Let $\gamma:\pi_1(N) \to G$ be an epimorphism onto a finite group $G$ and $\pi : N_{G} \to N$ the induced $G$-cover.
Then
\[ \Delta_{N}^{\gamma}=\Delta_{N_G,\pi_{*}}\in \Z[H]. \]
\end{lemma}

Finally, we need to rewrite the Alexander polynomial $\Delta_{N_G,\pi_{*}}$ in terms of the full
Alexander polynomial of $N_{G}$; their relation is the following.

\begin{proposition} \label{relalex} \cite{FV06}\cite{Tu01} Let $N$ be a 3--manifold with non--empty toroidal boundary,
and let $N_G$ be the $3$--manifold defined as above. Furthermore  let $\Delta_{N_G} \in
\Z[H_{G}]$ be the (ordinary multivariable) Alexander polynomial. Then we have the following
equality in $\Z[H]$:
\\ If $b_{1}(N_{G}) > 1$, then
\begin{equation}  \label{proone} \Delta_{N_G,\pi_*}   = \left\{ \begin{array}{ll} \pi_{*}(\Delta_{N_{G}})  & \mbox{if $b_{1}(N) >
     1$}, \\ \\ (t^{div{\pi_{*}}}-1)\pi_{*}(\Delta_{N_{G}}) & \mbox{if $b_{1}(N) = 1$,
$\im \, \pi_{*} = \langle t^{div\, \pi_{*}} \rangle$, $t \in H$ indivisible}.
\end{array} \right.
\end{equation} If $b_{1}(N_{G}) = 1$, then $b_{1}(N) = 1$ and
\begin{equation} \label{protwo} \Delta_{N_G,\pi_*} = \pi_{*} (\Delta_{N_{G}}). \end{equation}
\end{proposition}

%===================================================
\section{Proof of Theorems \ref{mainthm} and \ref{thm:fibered}}\label{section:proofs}

%===================================================
\subsection{Proof of Theorem \ref{mainthm}}
The topological ingredient of the proof is a result on the virtual Betti number of link exteriors.
This result can be deduced quite directly from \cite[Theorem~1.3]{CLR97}, but it is perhaps
appropriate, in order to illustrate its nature, to break down the proof to emphasize the role of
the JSJ decomposition of a link exterior.

We start with the following results.

\begin{theorem}\label{thm:lue}\cite{Lu88,Ko87}
Let $N$ be an irreducible $3$--manifold containing an essential torus or annulus $S$; up to a lift
to a finite cover, we can assume that $S$ is non--separating. Then $S$ is either the fiber of a
fibration over $S^1$, or the virtual Betti number $vb_1(N)$ of $N$ is infinite.
\end{theorem}

\begin{remark} Recall that having virtual Betti number $vb_1(N)$ infinite
 means that $N$ admits finite covers of arbitrarily large Betti number.
\textit{A priori}, the covers do not have to be regular: however, to any finite cover ${\hat N}$
with fundamental group $\hat{\pi}$ we can canonically associate a finite regular cover ${\bar
N}$ determined by the subgroup $\bar{\pi} := \cap_{p\in \pi_1(N)} p\hat{\pi}p^{-1}$. This
subgroup is clearly a normal subgroup of both $\hat{\pi}$ and $\pi_1(N)$. Also, since
$\hat{\pi}\subset \pi_1(N)$ is of finite index we see easily that $\bar{\pi}$ is in fact the
intersection of finitely many subgroups of $\pi_1(N)$ of finite index. Therefore
$\bar{\pi}\subset \hat{\pi} \subset \pi_1(N)$ is of finite index as well, and ${\bar N}$ is a
finite cover. From standard arguments, we have $b_1({\bar N})\geq b_1(\hat{N})\geq
b_1(N)$, so we can assume that $N$ admits finite \textit{regular} covers of arbitrarily large
Betti number. The set of left cosets $\pi_1(N)/{\bar \pi}$ is a finite group, that we
denote by $G$, hence ${\bar \pi}$ is the kernel of an epimorphism $\gamma : \pi_1(N) \to
G$, so that ${\bar N} = N_G$.  \end{remark}

\begin{theorem}\label{thm:clr1}\cite[Theorem~2.7]{CLR97}
Let $N$ be an irreducible $3$--manifold with non--empty incompressible boundary all of whose
components are tori. Suppose that the interior of $N$ has a complete hyperbolic structure of
finite volume. Then $vb_1(N) = \infty$. \end{theorem}

The topological ingredient in the proof of Theorem \ref{mainthm} is then the
following observation.

\begin{lemma}\label{prop:vb1}
Let $L=L_1\cup\dots\cup L_m\subset S^3$ be an oriented link which is
neither  the unknot nor the Hopf link. Then $vb_1(X(L))=\infty$.
\end{lemma}

\begin{proof}
First note that if $L$ is a split link, i.e. if $X(L)=S^3\sm \nu L$ is reducible, then
$\pi_1(X(L))$ maps onto a free group with two generators, which implies that $vb_1(X(L))=\infty$
(cf. e.g. \cite{Ko87}).

We can therefore now assume that $L$ is not a split link. In particular no component of $L$ bounds a disk in the
complement of the components. By Dehn's Lemma this implies that the
boundary of $X(L)$ is incompressible. As $X(L)$ is irreducible and
has boundary, $X(L)$ is Haken, hence it admits a geometric
decomposition along a (possibly empty) family of essential tori
$\mathcal{T}$. We will break the argument in subcases.

First assume that $\mathcal{T}$ is non-empty. Clearly $X(L)$ cannot be covered by a torus bundle
over $S^{1}$ since $X(L)$ has boundary. It therefore
follows from Theorem \ref{thm:lue} that $vb_{1}(X(L)) = \infty$.

Now assume that  $\mathcal{T}$ is empty. By Thurston's
geometrization of Haken manifolds we deduce that either $X(L)$ is
Seifert-fibered or the interior of $X(L)$ has a complete hyperbolic
structure of finite volume.

In the hyperbolic case, Theorem \ref{thm:clr1} asserts that
$vb_{1}(X(L)) = \infty$.

We are left with the Seifert--fibered case. The classification of Seifert links (see \cite[Chapter
II]{EN85}) shows that $L$ is the link obtained by removing $m$ fibers, regular or singular, from
the $(p,q)$--Seifert fibration of $S^{3}$, where $(p,q)$ are coprime integers or $(0,\pm 1)$.
Depending on the type of the orbifold quotient (see Jaco \cite[Chapter VIII]{Ja80}), $X(L)$ either
contains essential tori or is \textit{special}. In the former case, Theorem \ref{thm:lue} implies
$vb_{1}(X(L)) = \infty$ right away. If $X(L)$ is special, checking case--by--case, $L$ is either:
a (nontrivial) $(p,q)$--torus knot, obtained by removing a regular fiber; the union of the unknot
and its $(p,q)$--cable, obtained by removing a regular fiber and the fiber with multiplicity $p$
(whose exterior is the \textit{$p/q$--cable space}); one of a family of $3$--component links
obtained by removing a regular fiber and the two singular fibers. In the last two cases, we can
identify an essential, non--separating \textit{cabling} annulus joining a regular and a singular
fiber of the Seifert fibration. With the exception of the Hopf link with either orientations
(corresponding to $q = \pm 1$) these annuli do not fiber $X(L)$ by \cite[Theorem 11.2]{EN85}. For
a $(p,q)$--torus knot traced on a torus $T$, the annulus $X(L) \cap T$ is the only essential
annulus, and it is separating, so we pass to some finite cover. However, this cover cannot be an
annulus bundle over $S^1$ ($T^2 \times I$ or the twisted $I$-bundle over a Klein bottle), as by
\cite[Theorems 10.5, 10.6]{He76} the only manifolds covered by $T^2\times I$ are $T^2 \times I$
itself and the twisted $I$-bundle over a Klein bottle, which does not embed in $S^{3}$. It follows
that, with the exception of the Hopf link with either orientation (for whom $X(L) = T^2 \times
I$), all these links have $vb_{1}(X(L)) = \infty$ by Theorem \ref{thm:lue}.
\end{proof}

The following theorem, together with Proposition
\ref{prop:deltagdivides}, immediately implies  Theorem \ref{mainthm}.

\begin{theorem}\label{mainthm2}
Let $L=L_1\cup\dots\cup L_m\subset S^3$ be an oriented link which is
neither  the unknot nor the Hopf link. Then there exists an
epimorphism $\gamma:\pi_1(X(L))\to G$ onto a finite group $G$ such
that $\Delta^\gamma_L \ne \pm 1$.
\end{theorem}

\begin{proof}
Since $L$ is
neither  the unknot nor the Hopf link, Lemma \ref{prop:vb1} implies that there
exists a cover $X(L)_G$ with $b_{1}(X(L)_G) > 2$. As $X(L)_{G}$ has
non--empty boundary all of whose components are tori, Corollary II.4.4 of \cite{Tu02}
implies that the sum of the coefficients of
$\Delta_{X(L)_G}$ is zero. Hence, by Proposition \ref{relalex} the sum of
the coefficients of $\Delta_{X(L)_G,\pi_*}$ is zero as well, hence,
by Lemma \ref{lemmaalexg}, $\Delta_{X(L)}^{\gamma}$ cannot be $\pm
1$.
\end{proof}

When $L$ is the unknot or the Hopf link, $X(L)$
is homeomorphic to $S^1 \times D^2$ and $T^{2} \times I$ respectively.
In particular, the maximal (free) abelian cover ${\widehat X(L)}$
is contractible. Given any representation $\pi_1(X(L))\to \gl(R,k)$, we have
\[ H_1(X( L);R^k[H_1(X(L))]) \cong H_1({\widehat X(L)};R^k) = 0, \]
where the first isomorphism follows from the Eckmann--Shapiro lemma.
As the corresponding twisted Alexander module is trivial,
$\Delta_L^k=1$ for all $k$. This implies that the sequence
$\Delta_{L}^k$ detects the unknot and the Hopf link.

%===================================================
\subsection{Proof of Theorem  \ref{thm:fibered}}
Let $K \subset S^3$ be a fibered knot; it is shown in \cite{FK06} that $\Delta^\a_K$ is monic for
any representation $\a:\pi_1(X(K))\to \gl(\Z,k)$ (cf. also \cite{Ch03} and \cite{GKM05}). This
clearly implies that $\Delta^k_K$ is monic for all $k$.

Now let $K$ be a genus one knot such that
 $\Delta^k_K$ is monic for all $k$.
We denote by $N_K$ the zero framed surgery along $K$. Gabai \cite{Ga87} showed that $K$ is
fibered if and only if $N_K$ is fibered. Clearly, $N_K$ has vanishing Thurston norm. Under this hypothesis we show, in
\cite{FV06b} that if $N_K$ is not fibered, then there exists an epimorphism
$\beta:\pi_1(N_K)\to G$ onto a finite group $G$ such that $\Delta_{N_K}^\beta=0$.

Now consider the homomorphism
\[ \gamma: \pi_1(X(K))\to \pi_1(N_K)\to G.\]
Since $\pi_1(X(K))\to \pi_1(N_K)$ is an epimorphism, it follows from
the $5$--term exact sequence (cf. \cite[Chapter VII, Corollary
6.4]{B}) that $H_1(N_K;R^{k}[{t^\pm1}])$ is a quotient of
$H_1(X(K);R^k[t^{\pm1}])$, hence there exists a polynomial $p(t)\in
\zt$ such that $ \Delta_{K}^\gamma=p(t) \Delta_{N_K}^\beta$ (cf.
also \cite{KSW05}). In particular $ \Delta_{X(K)}^\gamma=0$. But
then Theorem  \ref{thm:fibered} follows from Proposition
\ref{prop:deltagdivides}.

%===========================================
\section{Calculations}\label{section:examples}

%===========================================
A natural test for invariants is their ability to detect mutation. Recall that two knots $K_1$ and
$K_2$ are called mutants if there exists a ball in $S^3$ whose boundary meets the knots in 4
points, such that removing the ball, rotating it by $\pi$ around an axis (in a way which preserves
the 4 points), and gluing it back turns $K_1$ into $K_2$.
\begin{figure}[h] \centering
\begin{tabular}{cc}
  \includegraphics[scale=0.45]{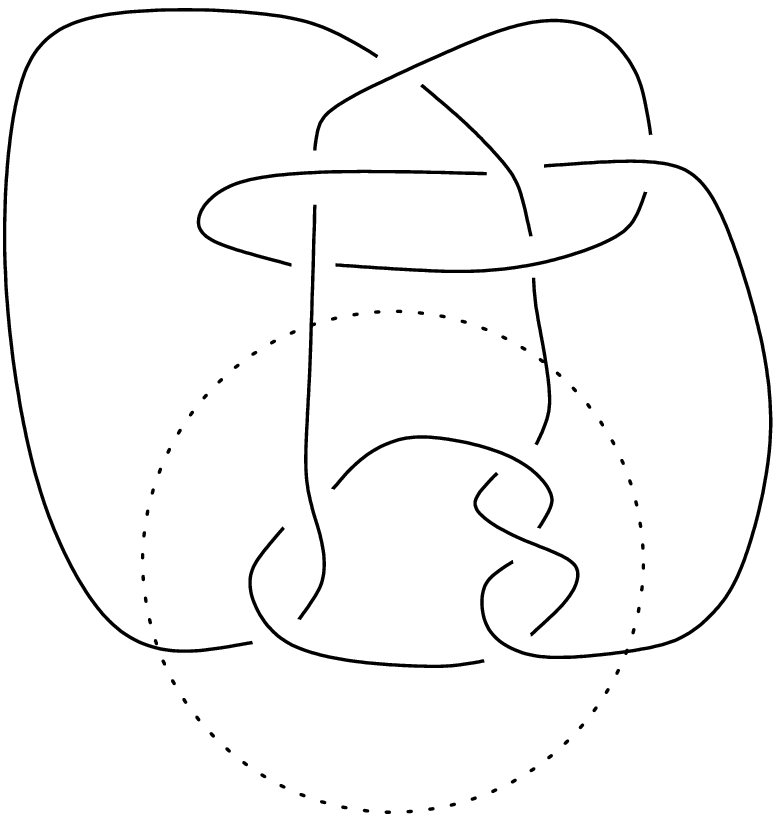}&
\includegraphics[scale=0.45]{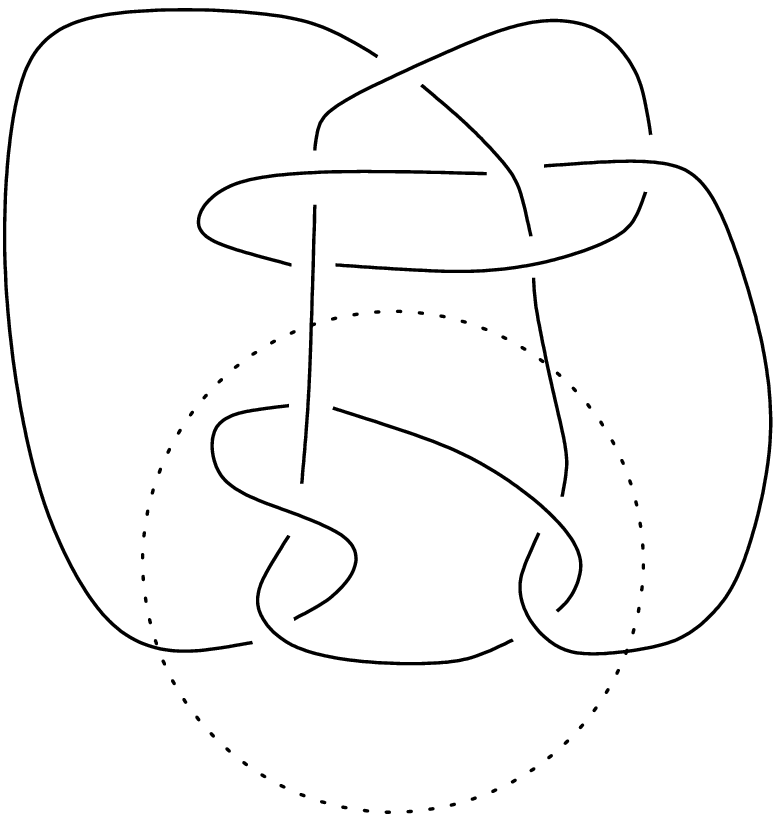}
\end{tabular}
\caption{The Conway knot and the Kinoshita--Terasaka knot.} \label{fig:mutant}
\end{figure}
Figure \ref{fig:mutant} shows perhaps the most famous pair of mutants, namely the  Conway knot
$11_{401}$ and the Kinoshita--Terasaka knot $11_{409}$. (Here we use Knotscape notation for knots
with more than 10 crossings, it is organized so that non-alternating knots are appended to
alternating ones instead of using `a' and `n' superscripts.) In both cases there exist seven
equivalence classes of abelian homomorphisms $\pi_1(X(K))\to S_5$ and one non--abelian equivalence
class of homomorphisms $\pi_1(X(K))\to S_5$. Using \textit{KnotTwister} we can compute their
invariants:
 \[ \ba{rcl}
\Delta_{11_{401}}^{5}\hspace{-0.15cm}&\hspace{-0.15cm}=\hspace{-0.15cm}&\hspace{-0.15cm}1+6t+9t^2+12t^3+t^5+3t^6+t^7+3t^8+t^9+12t^{11}+9t^{12}+6t^{13}+t^{14} \\
 \Delta_{11_{409}}^{5}\hspace{-0.15cm}&\hspace{-0.15cm}=\hspace{-0.15cm}&\hspace{-0.15cm}1+11t+12t^2+10t^3+5t^4+11t^5+4t^6+11t^7+5t^8+10t^9+12t^{10}+11t^{11}+t^{12}
 \ea
 \]
where both polynomials are considered in $\ftht$. Note
that Wada \cite{Wa94}  used parabolic representations to $\mbox{SL}(\F_{17},2)$ to show that these
two knots can be distinguished using twisted Alexander polynomials (cf. also \cite{In00}).

%It is well--known that invariants defined using a Skein relation formula do not detect mutants. In
%particular, the Alexander polynomial, the Conway polynomial, the Jones polynomial, the HOMFLY
%polynomial do not distinguish mutants. Furthermore, the computations of Shumakovitch using his
%program KhoHo  show that  any pair of mutant knots with 13 crossings or less has the same Khovanov
%homology. In fact no mutant pair of knots is known which can be distinguished by Khovanov homology.
%Other invariants which do not detect mutation are hyperbolic volume and the trace field of a
%hyperbolic knot. On the other hand knot Floer homology introduced by Ozsvath--Szabo \cite{OS04a}
%and Rasmussen \cite{Ra03} does distinguish (some) mutants. For example knot Floer homology
%determines the genus of a knot (cf. \cite{OS04b}) and the genus of the Conway knot equals three,
%whereas the genus of the Kinoshita--Terasaka knot equals two.

We have computed $\Delta^k_{K} \in \ftht$ for all groups of mutant
11--crossing knots for the smallest  value of $k$ that distinguishes
the mutants. The results are tabled in the Appendix. We also
computed $\deltaks5$ for all groups of mutant 12--crossing knots,
and again we verified that $\deltaks5$ distinguishes the mutant
knots.

We can summarize these computations in the following lemma.
\begin{lemma}
Let $K_1,K_2$ be a mutant pair with 12 crossings or less. Then
\[ \Delta^5_{K_1} \ne \Delta^5_{K_2}\in \F_{13}\tpm.\]
\end{lemma}

Note that the results of Section \ref{section:proofs} can be interpreted as stating that the
sequence $\Delta^k_K$ detects the unknot, the trefoil knot and the figure--8 knot (which are the
only fibered genus one knots). This raises the question about how effectively  the sequence
$\Delta^k_K$ at distinguishes knots in general.

In fact, we can use $\Delta^k_K$ to examine pairs of knots for whom other invariants are inconclusive.
For example,
the knots $10_{40}$ and ${10_{103}}$ are alternating knots with the same HOMFLY polynomial (hence
same Jones and Alexander polynomial) and the same signature. As Ng \cite[p.~292]{Ng05} points out
this implies by \cite{Le04} and \cite{OS03} that $10_{40}$ and ${10_{103}}$ also have the
 same Khovanov
homology and the same knot Floer homology. One can verify that
$\Delta_{10_{40}}^3=\Delta_{10_{103}}^3\in \F_{13}\tpm$ and that $\mathcal{R}_4(10_{40})$ and
$\mathcal{R}_4(10_{103})$ have eight elements each. Furthermore in $\F_{13}\tpm$ we have
\[ \ba{rcl}
\Delta_{10_{40}}^4&=&1+8t^2+t^3+12t^4+8t^5+\dots+8t^{176}+t^{178}\\
\Delta_{10_{103}}^4&=&1+11t+12t^2+4t^3+2t^4+3t^5+\dots+12t^{170}+11t^{171}+t^{172}.\ea
\] So the invariant $\Delta_K^4$ is neither determined by Khovanov
homology nor by Knot Floer homology.
%\footnote{Ng
%proposes the knots ${10_{56}})$ and ${10_{56}}$. They have $R_3=R_4=\emptyset$ and then have
%different number of elements in $R_5$, so one doesn't really need $\Delta^k_K$. In fact note, that
%the above two examples also have different $R_5$...} \footnote{
%  Ng
%\cite[p.292]{Ng05} says  that $11_{255}$ and $11_{257}$ are alternating 11--crossing knots with the
%same two--variable Kauffman polynomial. But at least according to Livingston's website they have
%different two--variable Kauffman polynomials.}
\newpage
\renewcommand{\theequation}{A-\arabic{equation}}
  % redefine the command that creates the equation no.
  \setcounter{equation}{0}  % reset counter
  \section*{Appendix}  % use *-form to suppress numbering

The following table lists all mutant pairs of knots with 11 crossings, together with the degrees
of $\Delta_K^{k}\in \ftht$ (for the smallest $k$ which distinguishes the mutants) and the first
5 terms of $\Delta_K^{k}$. All computations take place in the ring $\ftht$. Note that the first
five pairs of mutants are also distinguished by $\deltaks5$.
 \[
\ba{|l|l|l|l|l|l} \hline
\mbox{Knot} &k\phantom{t^{t^{t^t}}}&\#\mathcal{R}_k(K)&\deg(\Delta^k_K) & \mbox{Lowest and highest terms of }\Delta^k_K \in \ftht\\
\hline
 11_{44}&3&7&160& 1+12t+2t^2+7t^3+8t^4+12t^5+\dots+12t^{159}+t^{160}\phantom{t^{t^{t^t}}}\\
 11_{47}&3&7&160& 1+12t+11t^2+11t^3+7t^4+3t^5+\dots+12t^{159}+t^{160}\\
\hline
 11_{57}&3&7&160& 1+12t+5t^2+3t^3+2t^4+10t^5+\dots+12t^{159}+t^{160}\phantom{t^{t^{t^t}}}\\
 11_{231}&3&7&160&1+12t+t^2+7t^3+12t^4+7t^5+\dots+12t^{159}+t^{160}\\
\hline
11_{438}&3&7&118&1+11t+9t^2+4t^3+t^5+\dots+11t^{117}+t^{118} \phantom{t^{t^{t^t}}}\\
11_{442}&3&7&118&1+11t+3t^2+3t^3+t^4+10t^5+\dots+11t^{117}+t^{118}\\
\hline
 11_{440}&3&7&88&1+10t+t^2+12t^3+9t^4+5t^5+\dots+10t^{87}+t^{88}\phantom{t^{t^{t^t}}}\\
 11_{441}&3&7&76&1+10t+10t^2+11t^3+3t^4+2t^5+\dots+10t^{75}+t^{76}\\
 \hline
 11_{443}&3&7&160&1+2t+2t^2+8t^3+t^4+9t^5+\dots+2t^{159}+t^{160}\phantom{t^{t^{t^t}}}\\
 11_{445}&3&7&160&1+2t+6t^2+3t^3+6t^4+9t^5+\dots+2t^{159}+t^{160}\\
\hline
11_{19}&5\phantom{t^{t^{t^t}}}&13&496&1+9t+2t^2+12t^3+6t^4+6t^5+\dots+4t^{495}+12t^{496}\\
 11_{25}&5&12&460&1+2t+7t^2+9t^3+t^4+2t^5+\dots+2t^{459}+t^{460}\\
\hline
 11_{24}&5\phantom{t^{t^{t^t}}}&10&388&1+9t+11t^2+t^3+10t^4+t^5+\dots+9t^{387}+t^{388}\\
 11_{26}&5&9&352&1+2t+7t^2+8t^3+8t^4+8t^5+\dots+11t^{351}+12t^{352} \\
\hline
 11_{251}&5\phantom{t^{t^{t^t}}}&11&424&1+8t+7t^2+5t^3+12t^4+4t^5+\dots+5t^{423}+12t^{424}\\
 11_{253}&5&11&424&1+8t+3t^2+5t^3+t^4+9t^5+\dots+5t^{423}+12t^{424} \\
\hline
 11_{252}&5&9&352&1+2t+5t^3+6t^4+11t^5+\dots+11t^{351}+12t^{352}\phantom{t^{t^{t^t}}}\\
 11_{254}&5&10&388&1+9t+3t^2+2t^3+t^4+6t^5+\dots+9t^{387}+t^{388} \\
\hline
 11_{402}&5&17&466& 1+4t+3t^2+7t^3+8t^4+6t^5+\dots+9t^{465}+12t^{466}\phantom{t^{t^{t^t}}}\\
 11_{410}&5&15&418& 1+10t+10t^2+4t^3+2t^4+5t^5+\dots+3t^{417}+12t^{418}\\
\hline
 11_{403}&5&9&352&1+12t+4t^2+7t^3+7t^4+7t^5+\dots+t^{351}+12t^{352}\phantom{t^{t^{t^t}}}\\
 11_{411}&5&9&352&1+12t+3t^2+9t^3+2t^4+t^5+\dots+t^{351}+12t^{352}\\
\hline
 11_{406}&5&17&288&1+8t+2t^2+7t^3+2t^4+5t^5+\dots+8t^{287}+t^{288}\phantom{t^{t^{t^t}}}\\
 11_{412}&5&19&408&1+2t+t^2+12t^3+3t^4+5t^5+\dots+2t^{407}+t^{408}\\
\hline
 11_{407}&5&12&336& 1+7t+4t^2+5t^3+5t^4+6t^5+\dots+7t^{335}+t^{336}\phantom{t^{t^{t^t}}}\\
 11_{413}&5&12&340& 1+t+2t^2+t^4+3t^5+\dots+t^{339}+t^{340}\\
\hline
 11_{408}&5&15&568&1+4t+3t^3+6t^4+5t^5+\dots+4t^{567}+t^{568}\phantom{t^{t^{t^t}}}\\
 11_{414}&5&16&604&1+11t^2+10t^3+7t^4+5t^5+\dots+2t^{602}+12t^{604}\\
\hline
 11_{518}&5&12&220&1+t+12t^2+12t^3+3t^4+11t^5+\dots+t^{219}+t^{220}\phantom{t^{t^{t^t}}}\\
 11_{519}&5&11&228&1+3t+7t^2+11t^3+3t^4+5t^5+\dots+10t^{227}+12t^{228}\\
 \hline \ea \]

\end{document}